\documentclass{amsart}
\usepackage{pinlabel}
\usepackage{amsmath, amsthm, amssymb, amscd, mathrsfs, eucal, epsfig,color}
\usepackage{hyperref}

\begin{document}

\newtheorem{theorem}{Theorem}[section]
\newtheorem{lemma}[theorem]{Lemma}
\newtheorem{corollary}[theorem]{Corollary}
\newtheorem{conjecture}[theorem]{Conjecture}
\newtheorem{proposition}[theorem]{Proposition}
\newtheorem{question}[theorem]{Question}
\newtheorem*{answer}{Answer}
\newtheorem{problem}[theorem]{Problem}
\newtheorem*{claim}{Claim}
\newtheorem*{criterion}{Criterion}
\theoremstyle{definition}
\newtheorem{definition}[theorem]{Definition}
\newtheorem{construction}[theorem]{Construction}
\newtheorem{notation}[theorem]{Notation}
\newtheorem{convention}[theorem]{Convention}
\newtheorem*{warning}{Warning}
\newtheorem*{assumption}{Simplifying Assumptions}

\theoremstyle{remark}
\newtheorem{remark}[theorem]{Remark}
\newtheorem{example}[theorem]{Example}
\newtheorem{scholium}[theorem]{Scholium}
\newtheorem*{case}{Case}

\title{Dipoles and pixie dust}
\author{Danny Calegari}
\address{Department of Mathematics \\ University of Chicago \\
Chicago, Illinois, 60637}
\email{dannyc@math.uchicago.edu}
\date{\today}

\begin{abstract}
Every closed subset of the Riemann sphere can be approximated in the Hausdorff
topology by the Julia set of a rational map.
\end{abstract}

\maketitle

\section{Dipoles}
In \cite{Lindsey}, Kathryn Lindsey gives an elegant construction to show that any
Jordan curve in the complex plane can be approximated in the Hausdorff topology by
Julia sets of polynomials (another such construction was given subsequently by Oleg Ivrii
\cite{Ivrii}), and further, that any finite collection of disjoint
Jordan domains can be approximated by the basins of attraction of a rational map.
The proof depends on an interpolation result due to Curtiss \cite{Curtiss}. 

In this note, we give a direct geometric (and computation-free) proof that any closed
subset of the Riemann sphere can be approximated in the Hausdorff topology by the Julia
set of a rational map.

\medskip

In the theory of electomagnetism, a {\em dipole} refers to a pair of oppositely 
charged particles with charges of equal magnitude. The electric field of 
a dipole falls off at a rate of $1/r^3$ because of the approximate 
cancellation of the fields at distances large compared to the separation 
of the particles (unlike the usual inverse square law for a system with nonzero net
charge). The point is that the dipole is effectively electrically neutral on large
scales.

By analogy, we define a {\em dipole} to be a degree 1 rational function $(z-a)/(z-b)$
with a zero and pole at distinct non-zero $a$, $b$, for which $|a-b|$ is ``small''.
The dipole is uniformly close to 1 outside a small disk containing $a$ and $b$.
We now explain how to use dipoles to build designer Julia sets.

\medskip

Let's suppose we want to build a rational function whose Julia set approximates $X$,
a closed subset of the Riemann sphere. For simplicity, suppose $X$ is disjoint from the
unit circle. 

First, start with the map $f:z \to z^N$ where $N$ is some fixed big integer. The Julia
set of $f$ is the unit circle, and $f$ has $0$ and $\infty$ as superattracting fixed
points. 

Second, pick some finite collection $Y$ of discrete points (``pixels'')
which approximates $X$ closely
in the Hausdorff sense. We build a new rational function $g_\epsilon$ which is the product of $f$
with a dipole centered at each point in $Y$, where for each dipole the zero and pole are
within distance $\epsilon$ of each other. As $\epsilon \to 0$, the Julia sets of
$g_\epsilon$ converge uniformly in the Hausdorff topology to $\hat{Y}$, which is 
equal to the union of $Y$ together with its $N^k$th roots (i.e.\/
its preimages under $f$), together
with the unit circle. To see this, observe that the dynamics of $g_\epsilon$ converge
uniformly to $f$ on compact subsets of the complement of this set, while the presence of
the dipole near each point $y$ guarantees a point in the Julia set of $g_\epsilon$ near
$y$, and near its preimages under $f$. If $N$ is very big, $\hat{Y}$ is very close in
the Hausdorff topology to $Y \cup S^1$.

But now we are basically done: in place of $f:z \to z^N$ we could use a map
$f_{C,p}:(z-p) \to C(z-p)^N$ whose Julia set is an arbitrarily small circle centered at
an arbitrary point $p$ (for instance, a point in $X$). Repeating the construction above,
we get rational maps with Julia sets as close as we like in the Hausdorff sense
to $Y \cup p$ for $Y$ an arbitrary finite set and $p$ arbitrary.

\begin{figure}[htpb]
\centering
\includegraphics[height=1.5in]{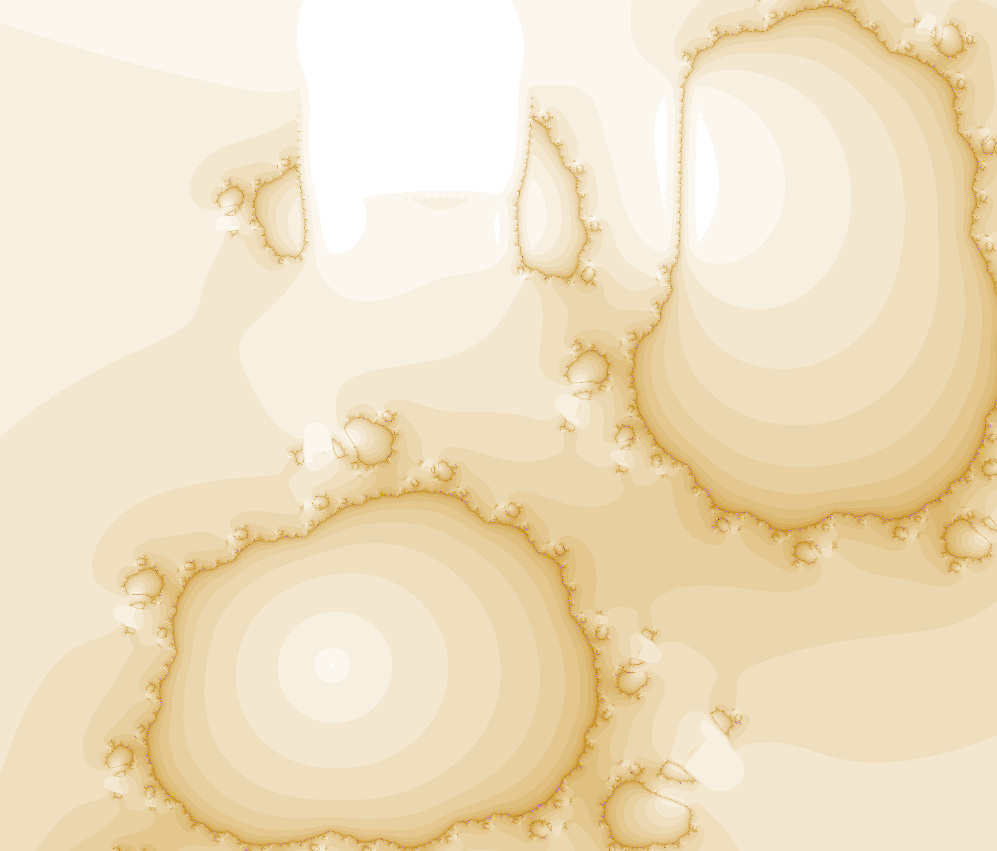}
\includegraphics[height=1.5in]{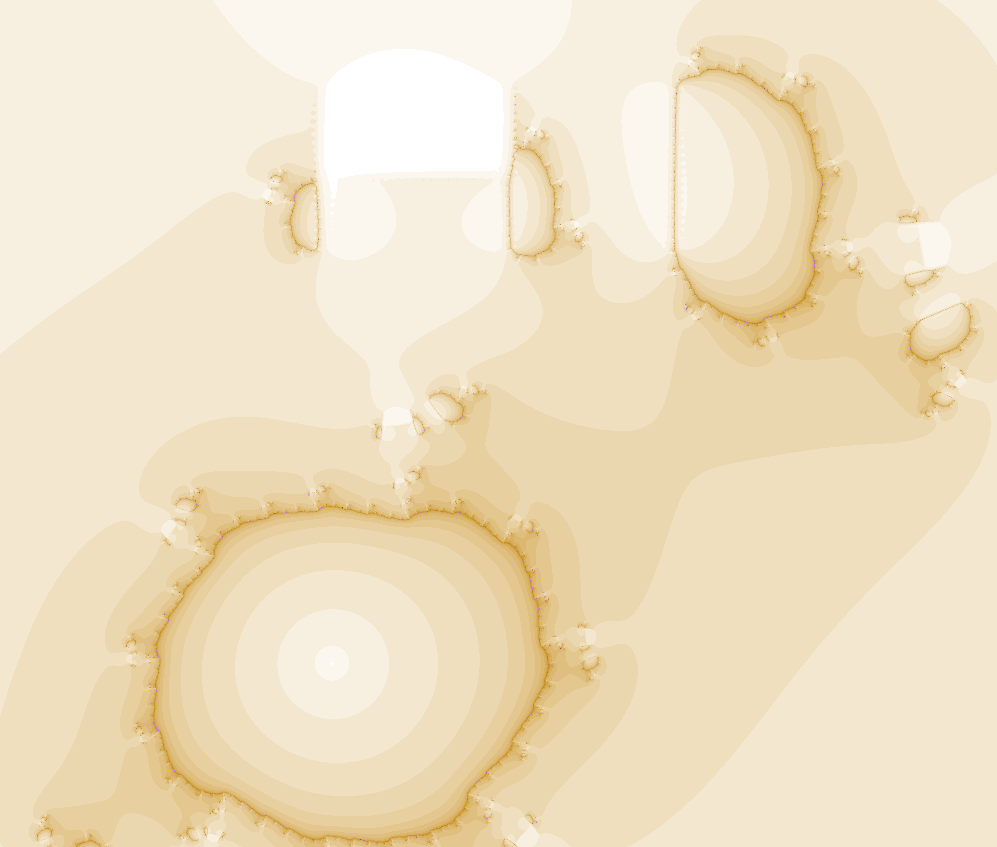}

\includegraphics[height=1.5in]{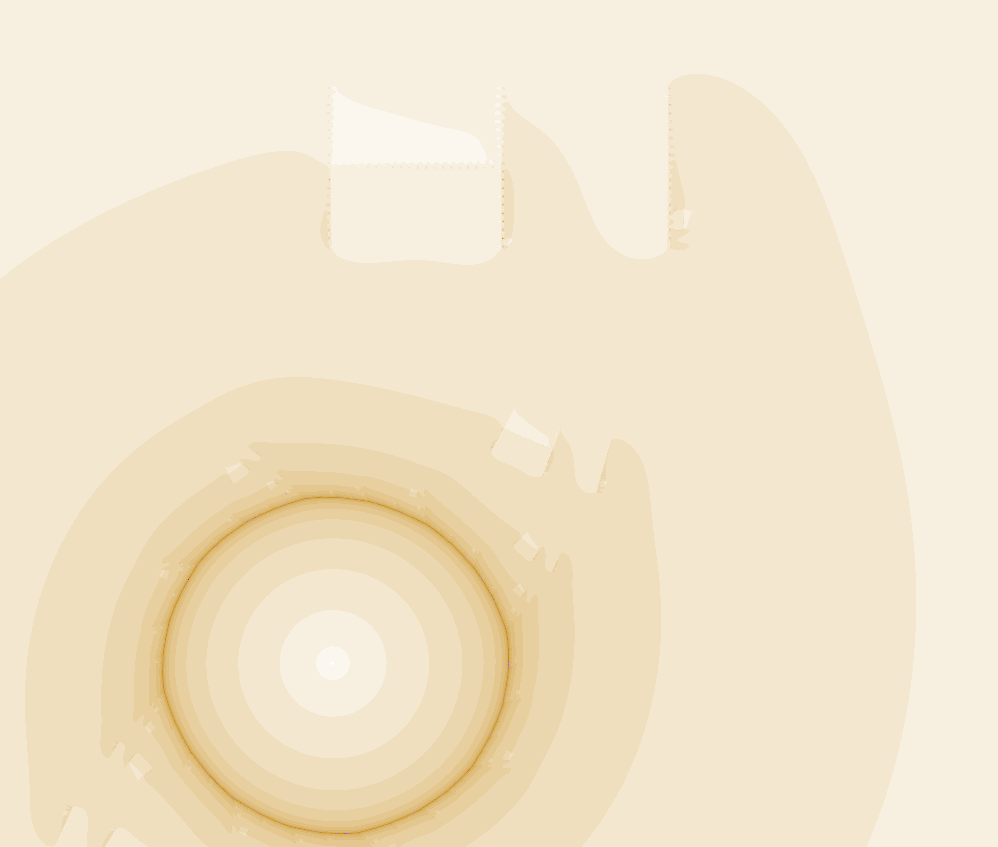}
\includegraphics[height=1.5in]{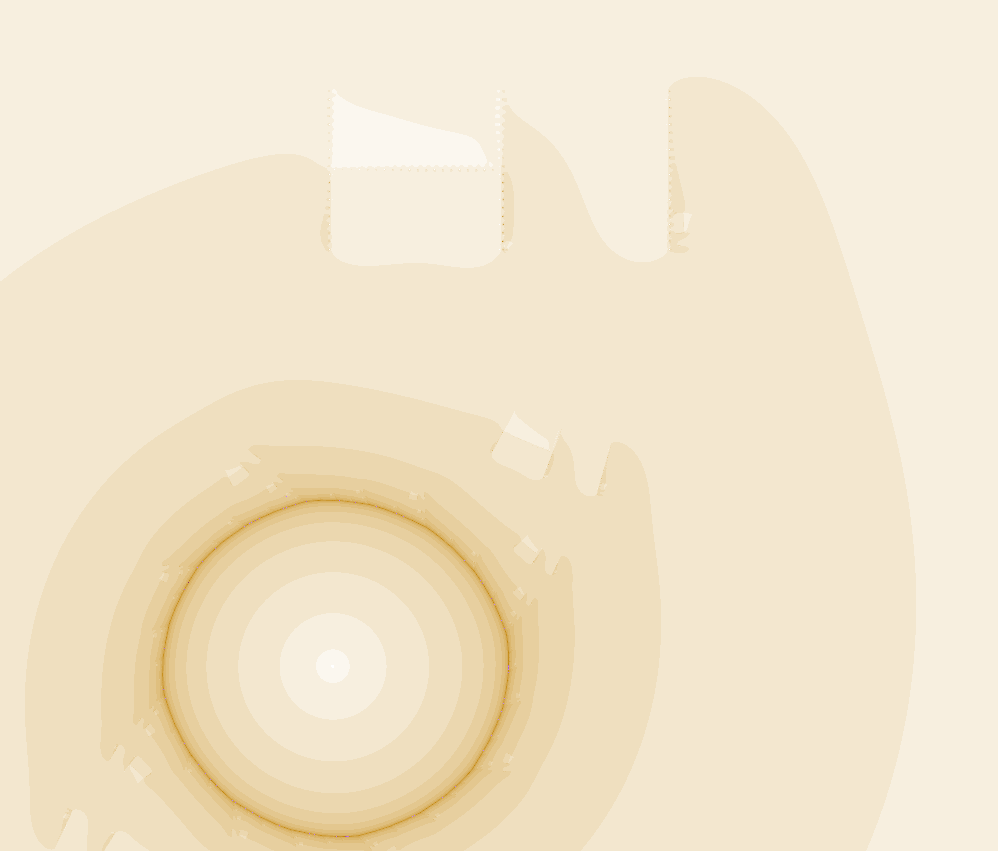}
\caption{80 dipoles and $N=2$ for $\epsilon = 0.2, 0.1, 0.05, 0.02$. The convergence
of the Julia sets to $\hat{Y}$ is evident. $Y$ is a pixelated ``HI'' at the top
of each picture.}\label{pixie}
\end{figure}

Figure~\ref{pixie} shows an example of a sequence of Julia sets 
constructed by this method.

\section{Acknowledgments}
I would like to thank Curt McMullen for bringing the reference \cite{Ivrii} to my attention.
Danny Calegari was supported by NSF grant DMS 1405466.

\end{document}